\documentclass{birkjour}
\usepackage{amssymb,amsmath,amsthm,enumerate,bbm}
\DeclareMathOperator{\supp}{supp}

\renewcommand\Im{\hbox{{\rm Im}}\,}

\numberwithin{equation}{section}

\theoremstyle{plain}
\newtheorem{theorem}{\bf Theorem}[section]
\newtheorem*{theorem*}{Theorem 1.1$'$}
\newtheorem{lemma}[theorem]{\bf Lemma}
\newtheorem{proposition}[theorem]{\bf Proposition}

\theoremstyle{definition}

\theoremstyle{remark}
\newtheorem*{remark*}{\bf Remark}
\newtheorem{remark}[theorem]{\bf Remark}

\newcommand{\wt}{\widetilde}

\DeclareFontFamily{U}{mathx}{\hyphenchar\font45}
\DeclareFontShape{U}{mathx}{m}{n}{<5> <6> <7> <8> <9> <10>
<10.95> <12> <14.4> <17.28> <20.74> <24.88> mathx10}{}
\DeclareSymbolFont{mathx}{U}{mathx}{m}{n}
\DeclareFontSubstitution{U}{mathx}{m}{n}
\DeclareMathAccent{\widecheck}{0}{mathx}{"71}
  \newcommand{\e}{\eqref}
\newcommand{\q}{\quad}

\newcommand{\ov}{\overline}

\newcommand{\1}{\mathbbm{1}}

\let\cal\mathcal

\let\Bbb\mathbb

\newenvironment{pf}{\begin{proof}}{\end{proof}}
\begin{document}

\title {Unbounded  Hankel  operators and moment problems}

\author{ D. R. Yafaev}
\address{ IRMAR, Universit\'{e} de Rennes I\\ Campus de
  Beaulieu, 35042 Rennes Cedex, FRANCE}
\email{yafaev@univ-rennes1.fr}

\begin{abstract} 
We   find  simple conditions for a non-negative Hankel quadratic form to be closable. Under some mild a priori assumption on the associated moments these sufficient conditions turn out to be also
   necessary. We also describe the domain of the corresponding closed form. This allows us to define unbounded non-negative Hankel operators under minimal assumptions on their matrix elements. The results obtained
  supplement the classical   Widom condition  for a Hankel operator to be bounded.
     \end{abstract}

\maketitle

%************************************************************
\section{Main results.  Discussion}  
%***********************************************************

{\bf 1.1.}
Hankel operators   $Q$ can formally be defined in the space $\ell^2({\Bbb Z}_{+})$ of sequences $g=(g_{0}, g_{1}, \ldots)$
 by the formula
  \begin{equation}
(Q g)_{n}=\sum_{m=0}^\infty  q_{n+m} g_m, \q n=0,1,\ldots.
 \label{eq:HD}\end{equation}
 Thus the matrix elements of a  Hankel operator depend  on the sum of indices only.
  
 The precise definition of the operator   $Q$ requires some accuracy.  Let  $\cal D \subset \ell^2({\Bbb Z}_{+})  $ be the dense set   of sequences $g=(g_{0},g_{1}, \ldots )$ with only a finite number of non-zero components. If the sequence 
 $q=(q_{0}, q_{1}, \ldots)\in \ell^2({\Bbb Z}_{+})  $, then for $g\in \cal D $ sequence \e{eq:HD} also belongs to $ \ell^2({\Bbb Z}_{+})  $. In this case the operator $Q$ is defined on $\cal D $, and it is symmetric  if $q_{n}= \ov{q_{n}}$. Without any a priori assumptions on $q_{n}$, 
  only the quadratic form
   \begin{equation}
q[g,g] =\sum_{n,m\geq 0} q_{n+m} g_{m}\bar{g}_{n}   
 \label{eq:QFq}\end{equation}
 is well defined for $g\in \cal D$.
 
 The fundamental theorem of Nehari \cite{Nehari} states that a Hankel operator $Q$ (defined, possibly, via its quadratic form \e{eq:QFq}) is bounded if and only if $q_{n}$ are the Fourier coefficients of some bounded function on the unit circle $\Bbb T$. The theory of Hankel operators is a very well developed subject. 
 We refer to the books \cite{NK, Pe} for basic information on this theory. However
 to the best of our knowledge, it was always assumed that Hankel operators were bounded. The only exception is paper \cite{Yf1a} where Hankel operators were realized as integral operators in the space $L^2({\Bbb R}_{+}) $.
 
 The goal of this paper is to make first steps in the study of {\it unbounded} Hankel operators. We consider {\it  non-negative} quadratic forms \e{eq:QFq} (in particular, we always assume that $q_{n}= \ov{q_{n}}$) so that
 we are tempted to define $Q$ as a self-adjoint operator corresponding to the quadratic form $q[g,g]$. Such an   operator exists if  the form $q[g,g]$ is closable in the space $\ell^2 ({\Bbb Z}_{+})$, but as is well known this is not always true.   We refer to the book \cite{BSbook} for  basic information concerning these notions.

 \medskip

{\bf 1.2.}
Below we give necessary and sufficient conditions guaranteeing the existence of a  closure of $q[g,g]$, but previously we discuss the link of Hankel quadratic forms with the Hamburger moment problem. The following result was obtained in \cite{Hamb}.
 
   \begin{theorem}\label{Hamb} 
   The condition 
      \begin{equation}
 \sum_{n,m\geq 0} q_{n+m} g_{m}\bar{g}_{n}\geq 0, \q \forall g\in \cal D,
 \label{eq:QF}\end{equation}
  is satisfied if and only if there exists a
  non-negative measure $dM(\mu)$ on $\Bbb R$ 
 obeying the condition
  \begin{equation}
 \int_{-\infty}^\infty | \mu|^n dM(\mu)<\infty, \q \forall n=0,1,\ldots,
 \label{eq:WH2}\end{equation}
such that the coefficients $q_{n }$ admit the representations 
   \begin{equation}
q_{n } = \int_{-\infty}^\infty \mu^n dM(\mu), \q \forall n=0,1,\ldots.
 \label{eq:WH}\end{equation}
    \end{theorem}

 Note that the measure satisfying equations \e{eq:WH} is in general not unique (see the paper \cite{Simon}, for a comprehensive discussion of this phenomenon). Roughly speaking,  the non-uniqueness of solutions of  the Hamburger moment problem is due to a very rapid growth of the  coefficients $q_{n}$.  Indeed, the famous Stieltjes example    shows that the measures
   \[
  dM_{\theta}(\mu)= \1_{{\Bbb R}_{+}} (\mu)\mu^{-\ln \mu} \big(1+\theta \sin (2\pi \ln \mu)\big) d\mu ,\q \forall \theta \in [-1,1],
 \]
 solve equations \e{eq:WH} with 
 $
 q_{n} =\sqrt{\pi}e^{(n+1)^2/ 4}.
$
 On the other hand, if
$
| q_{n}| \leq  R^n n!
$
 for some   $R>0$, then the solution of equations \e{eq:WH} for the measure $dM(\mu)$ is   unique.  
 
 \medskip

{\bf 1.3.}
The definition of the Hankel operator    requires essentially more restrictive assumptions which can be stated either in terms of the matrix elements $q_{n}$ or of the measure $dM(\mu)$.  

We say that a sequence $\varkappa_{n}>0$, $n\in {\Bbb Z}_{+}$, satisfies the Carleman condition if
\begin{equation}
\sum_{n\geq 1} \gamma_{n}^{-1} =\infty \q \mbox{where}\q \gamma_{n}=\inf_{m\geq n} m \varkappa_m.
\label{eq:Carl}\end{equation} 
Suppose that a function $f\in C^{\infty} ( \Delta)$ for some interval $\Delta\subset {\Bbb R} $ and
\[
|f^{(n)} (x)| \leq   n!  \varkappa_{n}^{n}, \q \forall n\in {\Bbb Z}_{+},
\]
where the sequence $ \varkappa_{n}$ obeys condition \e{eq:Carl}. Then (see the book \cite{Carleman}) $f$ is   quasi-analytic, that is,  the conditions $f^{(n)} (x_{0})=0$ for some $x_{0} \in \Delta $ and all $n\in {\Bbb Z}_{+} $ imply that 
$f(x)=0$ for  all $x\in \Delta$. 
If $ \varkappa_{n}= const$, then $f$ is an analytic function. The cases $ \varkappa_{n}= \varkappa_{0}\ln n$, $ \varkappa_{n}= \varkappa_{0}\ln n \ln (\ln n)$, etc., are known as the Denjoy conditions.

%We use below the Denjoy condition for the quasi-analyticity. Set $\ln_0 x= 1$, $\ln_{1} x= \ln x$ and  $\ln_{p+1} x= \ln (\ln_{p}x)$ for $p\geq 1$.  Then the estimates
%\[|f^{(n)} (x)| \leq A^{n} n! (\ln n \ln_2 n\cdots \ln_p n)^{n}\]
% for some $A\in {\Bbb R} $, some $p\geq 0$ and all $n\in{\Bbb Z}_{+}$ ensure that $f$ is quasi-analytic. The case $p=0$ leads to analytic functions.
%\begin{equation}
%| q_{n}| \leq A^{n} n! (\ln n \ln_{1} n\cdots \ln_p n)^{n}
%\label{eq:Denjoy}\end{equation} 

 Let us now state our main result.
 
   \begin{theorem}\label{Hamb1} 
     Let   assumption \e{eq:QF} be satisfied.  Then the following conditions are equivalent: 
     
     \begin{enumerate}[\rm(i)]
\item
The form $q[g,g]$ defined on $\cal D$ is closable in $\ell^2 ({\Bbb Z}_{+})$ and
\begin{equation}
q_{2n} \leq ( n! )^{2}  \varkappa_{n}^{2n}, \q \forall n\in {\Bbb Z}_{+},
\label{eq:Denjoy}\end{equation} 
for some sequence $ \varkappa_{n}$ obeying condition \e{eq:Carl}.

 \item
The matrix elements  $ q_n\to 0$ as $n\to \infty$.

\item
The measure $dM  (\mu)$ defined by equations  \e{eq:WH} satisfies  the condition
  \begin{equation}
M({\Bbb R}\setminus (-1,1) )=0
\label{eq:dense}\end{equation} 
$($to put it differently, $\supp M\subset [-1,1]$ and $M(\{-1\}) = M(\{1\})=0)$.
\end{enumerate}
        \end{theorem}
        
         \begin{remark}\label{rem} 

         \begin{enumerate}[\rm(i)]
         
           \item
           In the previous version of this article published as \cite{Y-mom}, condition
\e{eq:Denjoy} was omitted. It was pointed out in \cite{Berg} that, without some kind of an a priori assumption, the closability of $q[g,g]$ does not imply (ii)  or (iii).
         
         \item
Conditions \e{eq:Carl},
\e{eq:Denjoy} permit  very rapid growth of the moments $q_{n}$ as $n\to\infty$, for example, as $(n \ln n)^{n}$.

 \item
For $q\in \ell^2 ({\Bbb Z}_{+})$, the closability of $q[g,g]$ is obvious because in this case $q[g,g]=(Qg,g)$ where $Q$ is   the symmetric operator   defined on $\cal D$ by \e{eq:HD}.

\end{enumerate}
        \end{remark}

    As far as the proof of Theorem~\ref{Hamb1} is concerned, we note that only  the implication
        \[
        (i) \Longrightarrow (ii)\; {\rm or} \; (iii)
        \]
        is sufficiently non-trivial.     
        
            In Section~3 we  also give  (see Theorem~\ref{Clo}) an efficient description of the closure of the form \e{eq:QFq}.
             In Section~4 we discuss some consequences of our results for moment problems.
        
  \medskip

{\bf 1.4.}
          Theorem~\ref{Hamb1}  is to a large extent motivated by the following classical results of H.~Widom.
           
              \begin{theorem}\label{Widom}\cite[Theorem~3.1]{Widom}
                Let   the matrix elements  $ q_n $ of the Hankel operator \e{eq:HD}  be given by the relations 
                   \begin{equation}
q_{n } = \int_{-1}^1 \mu^n dM(\mu), \q \forall n=0,1,\ldots, \q M(\{-1\}) = M(\{1\})=0,
 \label{eq:WHq}\end{equation}
  with some non-negative  measure $dM  (\mu)$.
            Then the following conditions are equivalent: 
                
\begin{enumerate}[\rm(i)]
\item
The operator $Q$ is bounded.   
 \item
  $q_{n}=O(n^{-1})$ as $n\to \infty$.
   \item
 $M((1-\varepsilon,1))= O(\varepsilon)$ and  $M((-1, -1+\varepsilon))= O(\varepsilon)$ as $\varepsilon\to 0$.
\end{enumerate}
    \end{theorem}
    
     \begin{theorem}\label{Widom1}\cite[Theorem~3.2]{Widom}
     Under the same a priori assumptions as in    Theorem~\ref{Widom}
   the following conditions are equivalent: 
            
\begin{enumerate}[\rm(i)]
\item
The operator $Q$ is compact.   
 \item
  $ q_n=o (n^{-1})$ as $n\to \infty$.
   \item
 $M((1-\varepsilon,1))= o(\varepsilon)$ and  $M((-1, -1+\varepsilon))= o(\varepsilon)$ as $\varepsilon\to  0$.
\end{enumerate}
    \end{theorem}
    
    Theorems~\ref{Widom} 	and \ref{Widom1}
  give optimal conditions for the Hankel  operator $Q$ with matrix elements
 \e{eq:WHq}  to be bounded and compact.    Roughly speaking,    condition (iii) of   Theorem~\ref{Widom}    means that  the measure $dM(\mu)$ is ``subordinated" to the Lebesgue measure near the end points  $1$ and $-1$ of its support. Similarly,  condition (iii) of   Theorem~\ref{Widom1}
   means that  the measure $dM(\mu)$ 
  is  ``diluted"  compared to the Lebesgue measure near these end points.

  %************************************************************
  \section{Proof of Theorem~\ref{Hamb1}}  
%  
%* 

{\bf 2.1.}
 It is almost obvious that  conditions (ii) and (iii) are equivalent. Indeed, if 
(iii) is satisfied, then
 \[
q_{n } = \int_{-1+\varepsilon}^{1-\varepsilon} \mu^n dM(\mu)+ \int^{1 }_{1-\varepsilon} \mu^n dM(\mu)
+ \int_{-1 }^{-1+\varepsilon} \mu^n dM(\mu).
 \]
 The second and third integrals on the right are bounded by $M((1-\varepsilon,1))$ and $M((-1, -1 +\varepsilon))$, and hence they tend to zero as $\varepsilon\to 0$ uniformly in $n$. 
  The first integral is bounded by $ (1-\varepsilon )^n M((-1,1))$, and  therefore it tends to zero as $n\to \infty$ for every $\varepsilon>0$.
  Conversely, if there exists a set $X\subset {\Bbb R}\setminus (-1,1)$ such that $M(X)>0$, then $q_{2n}\geq M(X)$, and hence condition (ii) cannot be satisfied.

    It is convenient to reformulate   the fact that the form $q[g,g]$ is closable  in a different but equivalent form.  Let $ L^2 ( M)= L^2 (\Bbb R; dM)$ be the space of functions $u (\mu)$ with the norm $\| u\|_{ L^2 (M)}$.  We put
   \begin{equation}
({\cal A}g) (\mu)=\sum_{n=0}^\infty   g_{n} \mu^n
 \label{eq:A}\end{equation}
 and   observe that series \e{eq:A} converges for  all $g\in \ell^2 ({\Bbb Z}_{+})$ and all $\mu\in(-1,1)$. The function $({\cal A}g )(\mu)$ depends continuously on $\mu$, but only the estimate
 \[
 \sum_{n=0}^\infty |g_{n}| |\mu|^n \leq  (1-\mu^2)^{-1/2}\,  \| g\|_{\ell^2 ({\Bbb Z}_{+})} 
 \]
 holds. So it is of course possible that ${\cal A}g\not\in L^2 (M)$, but obviously 
 ${\cal A} g\in L^2 ( M)$ for all $g \in\cal D$.
Therefore we can  define   an auxiliary operator $A \colon \ell^2 ({\Bbb Z}_{+})\to   L^2 ( M)$  on domain $\cal D (A)=\cal D$ by the formula $Ag={\cal A}g$.

 In view of  equations  \e{eq:WH} the form    $   q[g,g]$  defined by relation  \e{eq:QFq}   can be written as
   \begin{multline}
   q[g,g] =\sum_{n,m=0}^\infty g_{m}\ov{g_{n}}  \int_{-\infty}^\infty  \mu^{n+m} dM(\mu)
   \\
    =\int_{-\infty}^\infty  | ({\cal A} g) (\mu)|^2 dM(\mu)= \| Ag\|^2_{L^2 ( M)}, \q g\in \cal D.
 \label{eq:A2}\end{multline}
 This yields the following result.
 
    \begin{lemma}\label{de}
 The form $   q[g,g]$ defined   on $\cal D$   is closable in the space $\ell^2 ({\Bbb Z}_{+})$  if and only if the operator $A \colon \ell^2 ({\Bbb Z}_{+})\to   L^2 ( M)$ defined   on the same set $\cal D$  
is closable.  
    \end{lemma}
    
      Recall that    the operator $A $   is closable if and only if its adjoint operator $A^* \colon L^2 ( M)\to \ell^2 ({\Bbb Z}_{+}) $ is densely defined. So our next goal is to construct $A^*$.  Observe that      under assumption \e{eq:WH2}
 for an arbitrary $u\in L^2 ( M)$, all the integrals
       \begin{equation}
 \int_{-\infty}^\infty  u(\mu) \mu^n dM(\mu)=: u_{n}, \q n\in {\Bbb Z}_{+},
 \label{eq:A1}\end{equation}
 are absolutely convergent. 
 We denote by   ${\cal D}_{*}\subset   L^2 (M)$ the set of   all $u \in   L^2 (M)$ such that the sequence $ \{u_{n}\}_{n=0}^\infty\in \ell^2 ({\Bbb Z}_{+})$.

   \begin{lemma}\label{LTM}
     Under assumption \e{eq:WH2} the operator $A ^*$ is given by the equality 
        \[
        (A^ {*}u)_{n}=\int_{-\infty}^\infty  u(\mu) \mu^n dM(\mu), \q n\in {\Bbb Z}_{+},
 \]
  on the domain ${\cal D}(A^*) ={\cal D}_{*}$. In particular, the operator $A $   is closable if and only if the set ${\cal D}_{*}$ is dense in $\ell^2 ({\Bbb Z}_{+})$.  
    \end{lemma}
    
     \begin{pf}
     Obviously, for all $g\in \cal D$ and all $u\in L^2 (M)$, we  have
       \begin{equation}
 (Ag,u)_{L^2 (M)} =   \sum_{n=0}^\infty g_{n}  \bar{u}_{n}
\label{eq:Y1}\end{equation} 
where the sequence $u_{n}$ is defined by relation \e{eq:A1}.
     The right-hand side here equals $(g, A^* u)$ provided $u\in {\cal D}_{*}$. It follows that $   {\cal D}_{*}\subset {\cal D}(A^*)  $.

         Conversely, if $u \in {\cal D}(A^*) $, then 
         \[
         |(Ag, u)_{L^2(M)} | =   |(g, A^* u)_{ \ell^2 ({\Bbb Z}_{+})} | \leq    \| g\|_{ \ell^2 ({\Bbb Z}_{+})}
         \, \| A^* u\|_{ \ell^2 ({\Bbb Z}_{+})} 
         \]
          for all $ g \in {\cal D}$. Therefore it follows from equality \e{eq:Y1} that
               \[
\big|   \sum_{n=0}^\infty g_{n}  \bar{u}_{n} \big| \leq C(u) \| g\|_{ \ell^2 ({\Bbb Z}_{+})}, \q \forall g \in {\cal D}.
\]
Since $ {\cal D}$ is dense in $ \ell^2 ({\Bbb Z}_{+})$, we  see that
 $\{ u_{n}\}_{n=0}^\infty \in \ell^2 ({\Bbb Z}_{+})$, and hence $u\in {\cal D}_{*}$.
          Thus   ${\cal D}(A^*)   \subset {\cal D}_{*}$.
          \end{pf}

 \medskip
 
 {\bf 2.2.}
Next, we use the following analytical  result.

     \begin{theorem}\label{dense}
     Under assumption \e{eq:dense}
     the set ${\cal D}_{*}$ is dense in the space $L^2 (M)$. Conversely,  if    the set ${\cal D}_{*}$ is dense in the space $L^2 (M)$ and  condition \e{eq:Denjoy} is true, then \e{eq:dense} is satisfied.
    \end{theorem}
    
      \begin{pf}
       Let the set $\cal E$ consist of $u\in L^2(M)$ such that $\supp u \subset [-a,a]$ for some $a<1$. According to definition 
       \e{eq:A1} for $u\in \cal E$, we have 
              $$
| u_{n}| \leq  a^n \| u \|_{L^2(M)} \, \sqrt{M ((-a,a))},
$$
whence ${\cal E}\subset{\cal D}_{*}$.
       Under assumption  \e{eq:dense} the set $\cal E$ is  dense in $L^2 (M )$ and so ${\cal D}_{*}$ is also dense in this space.

Let us prove the converse statement.   For an arbitrary $u\in L^2(M)$, we put
 \begin{equation}
f (x) = \int_{-\infty}^\infty e^{i \mu x} u(\mu) dM(\mu),\q x\in {\Bbb R}.
\label{eq:X1}\end{equation}
Then, for all $n\in{\Bbb Z}_{+}$, we have
 \begin{equation}
f^{(n)} (x) = i^{n}\int_{-\infty}^\infty e^{i \mu x} \mu^{n}u(\mu) dM(\mu)
\label{eq:XY}\end{equation}
and hence, by the Schwarz inequality,
 \begin{equation}
| f^{(n)} (x)  | \leq  \| u \|_{L^2(M)} \, \sqrt{ q_{{2n}} }.
\label{eq:X1a}\end{equation}
It now follows from condition \e{eq:Denjoy} that
  the function $f(x)$ is quasi-analytic on $\Bbb R$.

Assume now that $u\in{\cal D}_{*}$.
Then according to formula \e{eq:XY} for $x=0$ the sequence $f^{(n)} (0)$ is bounded and hence the function
\begin{equation}
\wt{f}(z):= \sum_{n=0}^\infty \frac{f^{(n)} (0)}{n!} z^n 
\label{eq:X4s}\end{equation}
is entire and  satisfies the estimate
\begin{equation}
|\wt{f}(z) |\leq C_{0}\sum_{n=0}^\infty \frac{1}{n!} |z|^n =C_{0} e^{|z|}, \q z\in {\Bbb C}, \q C_{0}=\max_{n\geq 0} |f^{(n)} (0)|.
\label{eq:X4}\end{equation}
Since $\wt{f}^{(n)}(0)=f^{(n)} (0)$ for all  $n\in{\Bbb Z}_{+}$ and both functions $\wt{f}(x)$ and $f(x)$  are quasi-analytic  on any bounded interval $\Delta\subset{\Bbb R}$,  they coincide for all $x\in{\Bbb R}$.

Let us now show that, for some $C>0$,
\begin{equation}
|\wt{f}(z)|\leq  C e^{| \Im {\textstyle z }|}, \q z\in {\Bbb C}.
\label{eq:X5}\end{equation}
Consider, for example the angle $\arg z\in [0,\pi/2]$ and put $F(z)= \wt{f}(z) e^{iz}$. Since $| e^{iz}|= e^{ -\Im {\textstyle z }}$,  it follows from  estimates \e{eq:X1a} and \e{eq:X4}  for $f(z)$ that $ |F(z)|\leq  C e^{ | z|}$ for all $z$ and that the function $F(z)$  is bounded on the rays $z=r$ and $z=ir$ where $r\geq 0$.
Therefore, by the Phragm\'en-Lindel\"of principle  (see, e.g., the book \cite{Evg}), $F(z)$  is bounded in the whole
angle $\arg z\in [0,\pi/2]$. This yields estimate \e{eq:X5} for $\wt{f}(z)= F(z)  e^{-iz}$.

 % Suppose that ${\cal D}_{*}$ is   dense in $L^2 (M )$.

According to the Paley-Wiener theorem (see, e.g., Theorem~IX.12 of \cite{RS}) it follows from estimate \e{eq:X5} that the Fourier transform of $\wt{f}(x)$ (considered as a distribution in the Schwartz class ${\cal S}' (\Bbb R)$) is supported by the interval $[ -1,1]$. Therefore  formula \e{eq:X1} for $f(x)=\wt{f}(x)$ implies that for every $u\in{\cal D}_{*}$, the distribution $u(\mu) dM(\mu)$ is also supported by  $[ -1,1]$,  that is
 \begin{equation}
 \int_{-\infty}^\infty \varphi(\mu)u(\mu) dM(\mu) =0, \q \forall \varphi\in C_{0}^\infty ({\Bbb R}\setminus [-1,1]).
\label{eq:X6}\end{equation}
If ${\cal D}_{*}$ is   dense in $L^2 (M)$, we can approximate $1$ by functions  $u\in {\cal D}_{*}$  in this space. 
Hence equality \e{eq:X6} is true with $u(\mu)=1$.   It follows that
 \begin{equation}
 \supp M \subset [-1,1]
\label{eq:X7}\end{equation}
because $\varphi\in C_{0}^\infty ({\Bbb R}\setminus [-1,1])$ is arbitrary.

For the proof of \e{eq:dense}, it remains to show that $M(\{-1\} ) =M(\{1\}) =0$. In view of \e{eq:X7} for an arbitrary $u\in L^2 (M )$, sequence \e{eq:A1}  admits the representation 
  \begin{equation}
    u _{n}=M(\{1\} ) u(1) + (-1)^n M(\{-1\} ) u(-1) + \int_{-1}^1 u(\mu) \mu^n d M_{0} (\mu) 
      \label{eq:X9}\end{equation}
 where $M_{0} (X)= M  (X\cap (-1,1))$ is the restriction of the measure $M$ on the open interval $(-1,1)$. Obviously, for any $\varepsilon\in (0,1)$, we have
     \begin{multline}
       \int_{-1}^1 u(\mu) \mu^n d M_{0} (\mu) =  \int_{-1+\varepsilon}^{1-\varepsilon} u(\mu) \mu^n d M_{0} (\mu)
\\  +
  \int_{1-\varepsilon}^{1 } u(\mu) \mu^n d M_{0} (\mu)  +
  \int_{-1 }^{-1+\varepsilon } u(\mu) \mu^n d M_{0} (\mu).
  \label{eq:X8}\end{multline}
  Applying the Schwarz inequality to each integral on the right, we estimate this expression  by
\[
\big((1-\varepsilon)^n   \sqrt{M_{0} ((-1,1))}+  \sqrt{M_{0} ((1-\varepsilon,1))}+ \sqrt{M_{0} ((-1 ,-1+\varepsilon))}\big) \| u\|_{L^2 (M_{0})} .
\]
Since $M_{0} ( (1-\varepsilon,1))\to 0$ and $M_{0}( (-1 ,-1+\varepsilon))\to 0$ as $\varepsilon\to 0$, we see that the integral 
in the left-hand side of \e{eq:X8} tends to zero as $n\to\infty$. Thus  \e{eq:X9} implies that
  \[
    u _{n}=M(\{1\} ) u(1) + (-1)^n M(\{-1\} ) u(-1) + o(1)
 \]
      as $n\to \infty$. 
 Therefore if $u\in  {\cal D}_{*}$, or equivalently $  \{ u_n\}_{n=0}^\infty \in \ell^2 ({\Bbb Z}_{+})$, then necessarily $M(\{1\} ) u(1) =M(\{-1\} ) u(-1) =0$. So if at least for one of the signs  $M(\{\pm 1\} )\neq 0$, then $  u(\pm 1) =0$. Let $u_{\pm}(\pm 1)=1$ and $u_{\pm}(\mu)=0$ for $\mu\neq \pm 1$. Since
 \[
 \| u_{\pm}- u\|_{L^2 (M)}\geq \sqrt{M(\{\pm 1\})}, \q \forall u\in  {\cal D}_{*},
 \]
  the function $u_{\pm}\in L^2 (M)$   cannot be approximated by functions    $ u\in {\cal D}_{*}$. Thus
 the set  $   {\cal D}_{*}  $ is not dense in $L^2 (M)$. 
   \end{pf}
 
% It follows that the function $v\in L^2 (M)$ such that $v(\mu)=0$ for $\mu\neq 1$ and $v(1)=1$ (or $v(\mu)=0$ for $\mu\neq -1$ and $v(-1)=1$) cannot be approximated by functions    $ u\in {\cal D}_{*}$.

        \begin{remark}\label{F-L}
        We have used   the Phragm\'en-Lindel\"of principle for the proof of estimate \e{eq:X5}  only. Actually, relation \e{eq:X6}
      can be directly deduced from estimates \e{eq:X1a} and \e{eq:X4} using the arguments given in the proof of Theorem~19.3 of the book \cite{Rudin}. However the intermediary estimate \e{eq:X5} makes the proof of 
 \e{eq:X6} essentially more transparent.
    \end{remark}

  \medskip
  
  {\bf 2.3.} 
    Let us come back to the proof of Theorem~\ref{Hamb1}. It remains to show that  the conditions (i) and (iii)  are equivalent.   In view of Lemmas~\ref{de} and \ref{LTM}  the form $q[g,g]$ is closable if and only if the set ${\cal D}_{*}$ is dense in $L^2 (M)$. Taking finally
  Theorem~\ref{dense} into account, we conclude the proof  of Theorem~\ref{Hamb1}.
        
         \begin{remark}\label{anal}
The condition \e{eq:Denjoy} in Theorem~\ref{Hamb1}  can be replaced by an estimate
\[
 \int_{-\infty}^\infty e^{2\epsilon|\mu|} dM(\mu) <\infty
\]
for some $\epsilon>0$. In this case   the function $f  (z)$ given by \e{eq:X1} is analytic and bounded in the strip  $|\Im z| <\epsilon$. Therefore the
functions $\wt{f}(z)$ (defined by \e{eq:X4s}) and $f(z)$ coincide as analytic   functions so that the theory of quasi-analytic functions is not required.
    \end{remark}

    \section{The closure of the Hankel quadratic form}  

{\bf 3.1.}
Let the condition \e{eq:dense} be satisfied. Recall that  the operator $\cal A$ was defined by equation  
 \e{eq:A} and $Ag={\cal A} g$ on domain ${\cal D}  (A) ={\cal D}$.
Let  $\bar{A}$ be  the closure of the operator $A$. According to Theorem~\ref{Hamb1} the form \e{eq:A2} is closable in the space $\ell^2 ({\Bbb Z}_{+})$ and the  form
 \begin{equation}
   q[g,g] = \| \bar{A} g\|^2_{ L^2 (M)}  
 \label{eq:AB2x}\end{equation}
 is closed  on domain ${\cal D}  [q] ={\cal D} (\bar{A})$. 
 
 Our goal is to find an efficient description  of ${\cal D}  [q]$. Since $\bar{A}=A^{**}$, we have to describe the set ${\cal D}(A^{**})$.   Let us define
   the operator $A_{\rm max}$   by the formula $A_{\rm max}g= {\cal A}g$   on the domain $\cal D(A_{\rm max})$ that consists of all $g\in \ell^2 ({\Bbb Z}_{+})$ such that ${\cal A}g\in L^2 (M)$. We will  show that
  \begin{equation}
A^{**} = A_{\rm max}.
 \label{eq:AB}\end{equation}

  A difficult part in the proof of \e{eq:AB} is the inclusion $A_{\rm max}  \subset A^{**} $ that is equivalent to the relation
          \begin{equation}
     (A_{\rm max}g,u)_{L^2 (M)} = (g, A^{*} u)_{\ell^2 ({\Bbb Z}_{+})}
 \label{eq:AB4}\end{equation}
 for all $g\in {\cal D} (A_{\rm max})$ and all $u\in {\cal D} (A^*)= {\cal D}_{*}$. In the detailed notation, relation \e{eq:AB4} means that
   \[
  \int_{-1}^1 \big(\sum_{n=0}^\infty g_{n}  \mu^n  \big) \ov{u(\mu)} dM(\mu)
  = \sum_{n=0}^\infty g_{n}\big(   \int_{-1}^1   \mu^n    \ov{u(\mu)} dM(\mu) \big). 
 \]
 The problem is that these series and integrals do not converge absolutely, and so the Fubini theorem cannot be applied.

\medskip

{\bf 3.2.}
     The shortest way to prove   \e{eq:AB} is to reduce
         the operator $\cal A$ by appropriate  unitary transformations   to the Laplace transform defined by the relation
   \begin{equation}
  ({\cal B} f) (\lambda)= \int_{0}^\infty e^{-t \lambda} f(t) dt.  
\label{eq:LAPj}\end{equation}
We consider it as a mapping ${\cal B} : L^2 ({\Bbb R}_{+}) \to L^2 ({\Bbb R}_{+}; d\Sigma)  $ where the non-negative measure $ d\Sigma (\lambda)$ on ${\Bbb R}_{+}$ satisfies the condition 
  \begin{equation}
\int_{0}^\infty (\lambda +1 )^{-k} d\Sigma (\lambda) <\infty 
 \label{eq:SS}\end{equation}
  for some $k>0$. The integral \e{eq:LAPj} converges for all $f\in L^2 ({\Bbb R}_{+})$ and $\lambda>0$. The function $ ({\cal B} f) (\lambda)$ is  continuous, but of course the estimate 
 \[
 |   ({ \cal B} f) (\lambda)|\leq (2\lambda)^{-1/2} \, \| f\|_{L^2 ({\Bbb R}_{+})}  
 \]
 does not guarantee that ${ \cal B}  f\in L^2 ({\Bbb R}_{+}; d\Sigma)$.

Let the set ${\sf D}\subset L^2 ({\Bbb R}_{+})$ consist of functions compactly supported in ${\Bbb R}_{+}$. If $f\in {\sf D}$, then $  ({\cal B} f) (\lambda)$ is a continuous function for all $\lambda\geq 0$ and $  ({\cal B} f) (\lambda)=O (e^{-c \lambda})$ with some $c=c(f)>0$ as $\lambda\to\infty$; in particular, ${\cal B} f\in L^2 ({\Bbb R}_{+})$. We put $Bf= {\cal B} f$ with ${\cal D} (B)= {\sf D}$. It is easy to show (see \cite{Yf1a}, for details) that the operator $B^*$ is given by the formula
 \begin{equation}
  (B^* v) (t)= \int_{0}^\infty e^{-t \lambda} v(\lambda) d\Sigma(\lambda) ,
\label{eq:Lad}\end{equation}
and its domain ${\cal D} (B^*)$ consists of all $v\in L^2 ({\Bbb R}_{+}; d\Sigma)$ such that $B^* v\in L^2 ({\Bbb R}_{+})$.
Obviously, this condition is satisfied if $v$ is compactly supported in ${\Bbb R}_{+}$.
Since the set of such $v$  is dense in $L^2 ({\Bbb R}_{+}; d\Sigma)$, the operator $B^*$ is densely defined. Thus $B$ admits the closure and $\bar{B}=B^{**}$.

  Let us now define
   the operator $B_{\rm max}$   by the formula $Bf= { \cal B} f$ on the domain $\cal D(B_{\rm max})$ that consists of all $f\in L^2 ({\Bbb R}_{+})$ such that ${ \cal B} f\in L^2 ({\Bbb R}_{+}; d\Sigma)$.  
   We use the following assertion.

\begin{lemma}\cite[Theorem~3.9]{Yf1a}\label{LLL}  
 Let     $ d\Sigma (\lambda)$ be a measure  on ${\Bbb R}_{+}$ such that the condition 
 \e{eq:SS} is satisfied for some $k>0$. Then
       \begin{equation}
{ B}^{**} = {  B}_{\rm max}.
 \label{eq:ABL}\end{equation}
        \end{lemma}

\medskip

{\bf 3.3.}
Let us find a relation between the operators $A$ and $B$.
Suppose that the   measures $ d\Sigma (\lambda)$ and $d M(\mu)$ are linked by the equality
  \begin{equation}
d M(\mu)=(\lambda +1/2)^{-2} d\Sigma (\lambda), \q \mu=\frac{2\lambda-1}{2\lambda+1}.
 \label{eq:MS}\end{equation}
 Thus if $M((-1,1))<\infty$, then the condition 
 \e{eq:SS} holds with $k=2$.
Put 
 \begin{equation}
  (Vu) (\lambda)= \frac{1}{\lambda+1/2} u \Big( \frac{2 \lambda-1}{2\lambda+1} \Big).  
\label{eq:LAV}\end{equation}
Obviously, $V: L^2 ((-1,1); dM) \to L^2 ({\Bbb R}_{+}; d\Sigma)$ is a unitary operator.
If the measures $dM $ and $d\Sigma$  are absolutely continuous, that is
   \[
  d\Sigma(\lambda)=  \sigma (\lambda) d\lambda, \q  \lambda>0, \q   dM(\mu)=\eta (\mu) d\mu , \q \mu \in (-1,1),
 \]
 with some $\sigma\in L^1_{\rm loc} ({\Bbb R}_{+})$ and $\eta\in L^1  (-1,1)$, then relation \e{eq:MS} means that
  \[
 \eta (\mu)=\sigma \big(\frac{ 1+\mu}{2 (1-\mu)}\big) .
 \]

  Recall that    the Laguerre polynomial (see the book \cite{BE}, Chapter~10.12) of degree $n$ is defined by the formula
\[
{\sf L}_{n}  (t)= n!^{-1} e^t   d^n (e^{-t} t^{n  }) / dt^n = \sum_{m=0}^n \frac{n! } {(n-m)! ( m!)^2} (-t)^m.
\]
We need  the identity (see formula (10.12.32) in \cite{BE})
 \begin{equation}
  \int_{0}^\infty     {\sf L}_{n} (t ) e^{-(1/2+ \lambda)t} dt=  \frac{1}{\lambda+1/2}\Big( \frac{2 \lambda-1}{2\lambda+1} \Big)^n, \q \lambda> -1/2 .
\label{eq:KL4x}\end{equation}
It can be deduced from this fact that the functions 
$  {\sf L}_{n}  (t) e^{-t/2}$, $n=0,1,\ldots$,   
form an orthonormal basis in the space $ L^2({\Bbb R}_{+})$, and hence the operator  $U \colon l^2({\Bbb Z}_{+})\to  L^2({\Bbb R}_{+})$ defined by the formula  
\begin{equation}
(U  g ) (t)=\sum_{n=0}^\infty g_{n}  {\sf L}_{n}  (t) e^{-t/2}, \q g=(g_{0}, g_{1}, \ldots),
\label{eq:K3}\end{equation}
is unitary.

A link between the operators ${\cal A}$ and ${\cal B}$ is stated in the following assertion. 

 \begin{lemma}\label{LA} 
     For all $g\in{\cal D}$, the identity  holds
      \begin{equation}
V {\cal A} g={\cal B} U  g.
\label{eq:LA}\end{equation}
        \end{lemma}
        
         \begin{pf}
It follows from equalities \e{eq:LAPj}, \e{eq:KL4x}   and \e{eq:K3}  that
   \[
( {\cal B}  U  g)(\lambda)  = \sum_{n=0}^\infty g_{n}   \int_{0}^\infty     {\sf L}_{n} (t ) e^{-(1/2+ \lambda)t} dt = \sum_{n=0}^\infty g_{n} \frac{1}{\lambda+1/2}\Big( \frac{2 \lambda-1}{2\lambda+1} \Big)^n.
 \]
In view of definitions  \e{eq:A}, \e{eq:LAV} the right-hand side here  equals $(V {\cal A}  g)(\lambda)$.
   \end{pf}
   
 Combining    Lemmas~\ref{LLL} and \ref{LA}, it is now easy to obtain
 the following result.
 
 \begin{lemma}\label{LM}  
 Let     $ d M (\mu)$ be a finite measure  on $(-1,1)$. Then    equality \e{eq:AB} holds.
        \end{lemma}
        
          \begin{pf}
          Observe that the adjoint of 
  the operator $B$ defined by   the formula $Bf= {\cal B}f$ on the set $  U{\cal D}$ is still given by formula  \e{eq:Lad}. Therefore it follows from \e{eq:LA} that $V   A =     B  U$, $   A^{*} V=  U   B^{*} $ and hence
        \begin{equation}
V   A^{**}=     B^{**} U.
\label{eq:LAx}\end{equation}

Let $g\in \ell^2 ({\Bbb Z}_{+})$ be arbitrary.
Approximating it by functions $g_{n}\in{\cal D}$ and using  \e{eq:LA}, we find  that $(V {\cal A} g)(\lambda)=({\cal B} U  g)(\lambda)$ for all $\lambda>0$. It follows that
 \begin{equation}
V A_{\rm max}  =  B_{\rm max} U   .
\label{eq:LAx1}\end{equation}
Comparing \e{eq:LAx} and \e{eq:LAx1}, we see that the identities \e{eq:AB} and \e{eq:ABL} are equivalent.
          \end{pf}

In view of  relation \e{eq:AB}, formula \e{eq:AB2x} leads to the following result.

  \begin{theorem}\label{Clo} 
     Let  the form $q[g,g]$ be defined on  ${\cal D}$ by  relation \e{eq:QFq}.
    Suppose that assumption \e{eq:QF} and one of three equivalent conditions (i), (ii) or (iii) of Theorem~\ref{Hamb1} are satisfied.       Let  ${\cal A}$ be  the operator  \e{eq:A}. 
    Then  the closure of $q[g,g]$
 is given by the equality 
 \begin{equation}
   q[g,g] =  \int_{-1}^1  | ( {\cal A}g) (\mu)|^2 dM(\mu) 
 \label{eq:AB2}\end{equation}
  on the set ${\cal D}[q]$ of all $g\in \ell^2 ({\Bbb Z}_{+})$ such that the right-hand side of \e{eq:AB2} is finite.
        \end{theorem}

        We recall that the non-negative operator $Q$ corresponding to the closed form \e{eq:AB2} satisfies the relations
        \begin{align*}
        q[g,h] &=  ( g, Q h), \q \forall g\in {\cal D}[q], \q \forall h\in  {\cal D} (Q)\subset {\cal D}[q],
        \\
        q[g,g]&= \| \sqrt{Q} g\|^2, \q \forall g\in  {\cal D}(\sqrt{Q}) = {\cal D}[q].
            \end{align*}
 Such an operator $Q$  is unique,           but its domain $ {\cal D} (Q) $ does not admit an efficient description.

            \section{Moment problems}

{\bf 4.1.}
Comparing  Theorems~\ref{Hamb} and \ref{Hamb1}, we obtain the following result concerning moment problems.

  \begin{proposition}\label{HambM} 
  A non-negative measure $dM(\mu)$ satisfying conditions  \e{eq:WH} and \e{eq:dense} exists  if and only if inequality 
  \e{eq:QF} holds and $q_{n}\to 0$ as $n\to\infty$ $($or, equivalently,   the form   \e{eq:QFq}  is closable and \e{eq:Denjoy}  is satisfied$)$.
    \end{proposition}

Instead of  the interval $[-1,1]$ we can consider an arbitrary finite  interval $[-a,a]$.  Our arguments proving the equivalence
of conditions (ii) and (iii) in  Theorem~\ref{Hamb1} lead to the following simple assertion.

  \begin{proposition}\label{HambM1} 
  A non-negative measure $dM(\mu)$ satisfying the condition 
     \[
q_{n } = \int_{-a}^a \mu^n dM(\mu), \q \forall n=0,1,\ldots,
 \]
  exists 
   if and only if inequality 
  \e{eq:QF} holds and $q_{n}= O (a^n)$ as $n\to\infty$. Moreover, $M(\{-a\})= M(\{a\})=0$ if and only if $q_{n}= o (a^n)$ as $n\to\infty$. 
    \end{proposition}
    
\medskip

{\bf 4.2.}
    The results obtained above can be combined with the Stieltjes theorem which states that   there exists a non-negative measure $dM(\mu)$ satisfying equations  \e{eq:WH} and 
 such that
$\supp M \subset [0,\infty )$
if and only if inequalities 
  \e{eq:QF} and
   \begin{equation}
 \sum_{n,m\geq 0} q_{n+m+1} g_{m}\bar{g}_{n}\geq 0, \q \forall g\in \cal D,
 \label{eq:QFst}\end{equation}
 hold.
 
 Let us state   analogues of Propositions~\ref{HambM}  and \ref{HambM1}. 
 
  \begin{proposition}\label{HambS} 
  A non-negative measure $dM(\mu)$ satisfying conditions  \e{eq:WH} and $\supp M \subset [0, 1]$,  $  M(\{1\})=0$
 exists  if and only if inequalities 
  \e{eq:QF} and   \e{eq:QFst}  hold and $q_{n}\to 0$ as $n\to\infty$ $($or, equivalently,   the form   \e{eq:QFq}  is closable and \e{eq:Denjoy}  is satisfied$)$.
    \end{proposition}
 
   \begin{proposition}\label{HambMst} 
  A non-negative measure $dM(\mu)$ satisfying the condition
     \[
q_{n } = \int_0^a \mu^n dM(\mu), \q \forall n=0,1,\ldots.
 \]
exists   if and only if inequalities 
  \e{eq:QF} and   \e{eq:QFst}  hold and $q_{n}= O (a^n)$ as $n\to\infty$. Moreover, $  M(\{a\})=0$ if and only if $q_{n}= o (a^n)$ as $n\to\infty$. 
    \end{proposition}

Note that the moment problem  \e{eq:WH} with the measure $dM(\mu)$ supported by  a compact interval is  called the Riesz problem.  The necessary and sufficient conditions for the existence of its solution    are well known (see, e.g., the book \cite{AKH}),
but they are stated in quite different terms compared to the results of this section.

 %%%%%%%%%%%%%%%%%%%%%%%%%%%%%%%%%%%%%%%%
%%%%%%%%%%%%%%%%%%%%%%%%%%%%%%%%%%%%%%%%

 \end{document}